\documentclass[ijoc,nonblindrev]{informs3} 

\OneAndAHalfSpacedXII 

\usepackage{natbib}
 \bibpunct[, ]{(}{)}{,}{a}{}{,}%
 \def\bibfont{\small}%
\TheoremsNumberedThrough     

\EquationsNumberedThrough    

\usepackage{rotating}
\usepackage[theme=grayscale]{jlcode}
\usepackage[utf8]{inputenc}

\usepackage{hyperref}
\usepackage[anythingbreaks]{breakurl}

\usepackage{array,multirow,graphicx}
\usepackage{float}
\usepackage{tabularx,dcolumn,booktabs}

\usepackage{subcaption}

\usepackage{tikz}

\DeclareMathOperator{\st}{s.t.}

\begin{document}


\RUNAUTHOR{Dias Garcia, Bodin and Street}

\RUNTITLE{Comparing \texttt{BilevelJuMP.jl} Formulations: Support Vector Regression Hyperparameter Tuning}

\TITLE{Comparing \texttt{BilevelJuMP.jl} Formulations:\\ Support Vector Regression Hyperparameter Tuning}

\ARTICLEAUTHORS{%
\AUTHOR{Joaquim Dias Garcia}
\AFF{LAMPS at PUC-Rio \& PSR, \EMAIL{joaquim@psr-inc.com}}
\AUTHOR{Guilherme Bodin}
\AFF{LAMPS at PUC-Rio \& PSR, \EMAIL{guilhermebodin@psr-inc.com}}
\AUTHOR{Alexandre Street}
\AFF{LAMPS at PUC-Rio, \EMAIL{street@ele.puc-rio.br}}
} 

\ABSTRACT{%
In this technical report, we compare multiple reformulation techniques and solvers that can be used with the Julia package BilevelJuMP. We focus on the special case of Hyperparameter Tuning  for Support Vector Regression. We describe a bilevel model for the problem in question. Then we present code for generating data and models that solve the problem. Finally, we present results and a brief analysis.
}%



\maketitle

%


\section{Introduction}\label{intro}

In this technical report, we compare multiple reformulation techniques and solvers that can be used with the Julia package \texttt{BilevelJuMP.jl} \cite{diasgarcia2022bileveljump}. We focus on the special case of Hyperparameter Tuning  for Support Vector Regression (HT-SVR).

We do not recommend the usage of the models presented here for actual HT-SVR. Although the formulation is precise, the solution method is not practical nor scalable. This problem is used as an example since it is easy to generate data for and vary the number of variables and constraints easily.

In Section \ref{intro}, we describe a bilevel model for the problem in question. Then, in Section \ref{code}, we present code for generating data and models that solve the problem. Finally, we present results in Section \ref{exp} and a brief analysis in Section \ref{ana}.

\section{Bilevel Optimization Formulation for HT-SVR}\label{intro}

In this section, we describe an interesting example of a bilevel program. The main goal is to start from a non-trivial problem, model it in \texttt{BilevelJuMP.jl}, and solve it with multiple methods to have a glimpse of the difference between techniques and solvers.

Hyperparameter tuning with bilevel optimization is a recent trend in the intersection of the Machine Learning and Optimization communities \citep{franceschi2018bilevel, kunisch2013bilevel, mackay2019self}. Although most hyperparameter tuning methods are based on bilevel optimization, the state-of-the-art solution methods are usually heuristic with special considerations to the problem in question. This problem is a good case to compare techniques due to the simplicity of the model and because small enough instances can be solved by standard bilevel optimization methods.

The example will follow the one from \cite{bennett2006model}, though with some simplifications. Given two data sets $O$ and $I$ with out-of-sample and in-sample data represented by the points labeled by $i$: $(y_i, \{x_{ij}\}_{j \in J})$, where $J$ is the set of features.
\[
\begin{aligned}
    &\min_{C \geq 0, \varepsilon \geq 0, \xi^U \geq 0} && \sum_{i \in O} \xi_i^U \\
    &\st && \xi_i^U \geq  + y_i - \sum_j w_j x_{ij}, \quad i \in O \\
    &    && \xi_i^U \geq - y_i + \sum_j w_j x_{ij} , \quad i \in O \\
        &&& w(C, \varepsilon) \in
     \begin{aligned}[t]
        &\argmin_{\xi^L \geq 0, w } && ||w||^2_2 + C \sum_{i \in I} \xi_i^L\\
            &\st && \xi_i^L + \varepsilon \geq  + y_i - \sum_{j \in J} w_j x_{ij} , \quad i \in I\\
            &    && \xi_i^L + \varepsilon \geq - y_i + \sum_{j \in J} w_j x_{ij}, \quad i \in I\\
     \end{aligned}\\
\end{aligned}
\]
The lower-level model is responsible for obtaining the best possible support vectors, $w$, given the problem data and the hyperparameters $C$ and $\varepsilon$. The two hyperparameters are variables selected by the upper level so that the support vector, $w$, optimized by the lower level has a minimal out-of-sample error.
The variables $\xi^U$ and $\xi^L$ denote the absolute value loss in the upper and lower models, respectively. The upper level is a linear program, while the lower level is a quadratic program.

\section{\texttt{BilevelJuMP.jl} Code for HT-SVR}\label{code}

In Figure \ref{fig:jump_svr}, we present \texttt{BilevelJuMP.jl} code to model the hyperparameter tuning of SVR described above. Thanks to the JuMP syntax, the code greatly resembles the abstract model, simplifying the writing and documenting of the code. 
\begin{figure*}[!ht]
    \centering
\begin{lstlisting}[language = Julia]
using JuMP, BilevelJuMP
# sample data
Features = 2
Samples = 10
J = 1:Features
I = 1:div(Samples, 2)
O = (div(Samples, 2)+1):Samples
x = 2 * (rand(Samples, Features) .- 0.5)
w_real = ones(Features)
y = x * w_real .+ 0.1 * 2 * (rand(Samples) .- 0.5)
# model building
model = BilevelModel()
@variable(Upper(model), C >= 0)
@variable(Upper(model), eps >= 0)
@variable(Upper(model), xi_U[i=O] >= 0)
@variable(Lower(model), w[j=J])
@variable(Lower(model), xi_L[i=I] >= 0)
@objective(Upper(model),
    Min, sum(xi_U[i] for i in O))
@constraints(Upper(model), begin
    [i in O], xi_U[i] >= + y[i] - sum(w[j]*x[i,j] for j in J)
    [i in O], xi_U[i] >= - y[i] + sum(w[j]*x[i,j] for j in J)
end)
@objective(Lower(model),
    Min, sum(w[j]^2 for j in J) + C * sum(xi_L[i] for i in I))
@constraints(Lower(model), begin
    [i in I], xi_L[i] + eps >= + y[i] - sum(w[j]*x[i,j] for j in J)
    [i in I], xi_L[i] + eps >= - y[i] + sum(w[j]*x[i,j] for j in J)
end)
\end{lstlisting}
    \caption{Code to model SVR hyper-parameter tuning}
    \label{fig:jump_svr}
\end{figure*}

That same code was used to perform a series of comparisons between solvers. We started by creating instances with a different number of features and observations (dataset size). We randomly created the matrix $x$ with a uniform distribution in $[-1, +1]$, then we created the \textit{real} $w$ as a vector of ones with appropriate dimension. Next, we defined $y = x w + \epsilon$, where $\epsilon$ follows a uniform distribution in $[-0.1, +0.1]$. Half of the dataset was considered in-sample data, while the other half was considered out-of-sample data. It is not our intention to be fully realistic here, our goal is to provide a didactic example.

\section{Experiments}\label{exp}

We created instances with $10$, $100$ and $1000$ samples. For all these numbers of samples, we created datasets with $1$, $2$ and $5$ features. For the cases with $100$ samples, we also created datasets with $10$, $20$ and $50$ features.

Finally, we optimized the bilevel problem for each data set with multiple reformulations and multiple solvers. The only solver attribute we set was a time limit of $600$ seconds ($10$ minutes) and left all other attributes as default, which might differ considerably from one solver to the other. 

We used
Julia 1.6.2,
BilevelJuMP 0.5.0,
CPLEX 22.1 \citep{cplex},
Gurobi 9.5 \citep{Gurobi},
HiGHS 1.2 \citep{huangfu_2018},
Ipopt 3.14 \citep{wachter2006implementation},
Knitro 13.0 \citep{nocedal2006knitro},
SCIP 8.0 \citep{BestuzhevaEtal2021OO},
Xpress 8.13 \citep{xpress}. All the required code, including exact package versions (see the \textit{manifest.toml} file),  is in the benchmarks folder of the git companion repository:
\url{https://github.com/joaquimg/BilevelJuMPBenchmarkSVRHT}.

We present the results in tables with a similar format. The first column describes the instance, the first number being the sample size and the second the number of features. Then we have three columns for each solver, the first, \textit{Obj}, presents the upper-level objective value returned by the solver (typically the best incumbent solution), the second contains the \textit{Gap} in percent (\%), Ipopt and KNITRO will not have gaps as they are NLP solvers, if no gap was reported the entry will be blank (with a ``$-$''), the third is \textit{Time} in seconds, if the time reaches $600$ the entry will be blank (with a ``$-$'').

Multiple \texttt{BilevelJuMP.jl} modes were considered as different techniques to solve bilevel problems. For more information on each of the specific modes, the reader is referred to the preprint \cite{diasgarcia2022bileveljump} and the user manual: \url{https://joaquimg.github.io/BilevelJuMP.jl/v0.6.1/tutorials/modes/}.
Table \ref{table_sos1} presents results for  \texttt{SOS1Mode} and \texttt{IndicatorMode}. Table \ref{table_fa100} presents results for \texttt{FortunyAmatMcCarlMode}, with big-Ms set to $100$, \texttt{StrongDualityMode} and \texttt{ProductMode}, the latter two with binary expansions, so the resulting problem is a MIP. The variable bounds were set to $+/-100$ for binary expansions. Finally, Table \ref{table_prod} presents the solutions of both \texttt{ProductMode} and \texttt{StrongDualityMode} for Non-Linear Programming solvers and Gurobi with its \textit{NonConvex} mode activated.

\section{Brief Analysis of Results}\label{ana}

We can make some comments and remarks based on the results in the tables. We note that \texttt{SOS1Mode} and \texttt{IndicatorMode} perform well in smaller instances, with a slight advantage for \texttt{SOS1Mode}. Interestingly, CPLEX's solution for $1000/01$ with \texttt{IndicatorMode} slightly disagrees with all solutions from other solvers with the \texttt{SOS1Mode}. \texttt{FortunyAmatMcCarlMode} and \texttt{StrongDualityMode} seem very amenable to MIP solvers. In particular, Gurobi closed the gap within the given 10 minutes for all but one instance in the latter mode. However, we must be careful since we selected arbitrary bounds for those methods. Moreover, \texttt{StrongDualityMode} also relies on binary expansion approximations, which led solvers to a solution that disagrees with the other methods on the $10/05$ instance. On the other hand, \texttt{ProductMode} is the worst strategy for MIP solvers in these instances. For NLP solvers, both \texttt{ProductMode} and \texttt{StrongDualityMode} return objective values that are close to the ones found by MIP solvers, but, in this case, there is a slight advantage for \texttt{ProductMode}. Finally, Gurobi NonConvex seems to work much better with \texttt{StrongDualityMode}, claiming very good results in the instances with $1000$ samples that agree with some of the other presented objective values.

The results are particular to a toy problem. However, the tables demonstrate that it is possible to benchmark multiple solvers and methods. Moreover, such benchmarks are easy to perform thanks to \texttt{BilevelJuMP.jl}.

\begin{table}[!ht]
\centering
\resizebox{13cm}{!}{
\begin{tabular}{rr|rrr|rrr|rrr|rrr}
\toprule
 &  & \multicolumn{3}{c}{CPLEX} \vline & \multicolumn{3}{c}{Gurobi} \vline & \multicolumn{3}{c}{SCIP} \vline & \multicolumn{3}{c}{Xpress}  \\
 & Inst & Obj & Gap  & Time  & Obj & Gap  & Time  & Obj & Gap  & Time  & Obj & Gap  & Time  \\
\midrule
 & 10/01  & 0.30 & 0 & 0 & 0.30 & 0 & 0 & 0.30 & 0 & 0 & 0.30 & 0 & 0 \\
 & 10/02  & 0.22 & 0 & 0 & 0.22 & 0 & 0 & 0.22 & 0 & 0 & 0.22 & 0 & 0 \\
 \parbox[t]{2mm}{\multirow{3}{*}{\rotatebox[origin=c]{90}{\texttt{SOS1Mode}}}}
 & 10/05  & 0.09 & 0 & 0 & 0.09 & 0 & 0 & 0.09 & 0 & 0 & 0.09 & 0 & 0 \\
 & 100/01  & 2.42 & 0 & 0 & 2.42 & 0 & 0 & 2.42 & 0 & 0 & 2.42 & 0 & 0 \\
 & 100/02  & 2.40 & 4 &  -  & 2.40 & 4 &  -  & 2.40 & 4 &  -  & 2.40 & 4 &  -  \\
 & 100/05  & 2.30 & 6 &  -  & 2.30 & 6 &  -  & 2.31 & 6 &  -  & 54.87 &  -  &  -  \\
 & 100/10  & 8.54 & 392 &  -  & 79.59 &  -  &  -  & 79.59 &  -  &  -  & 79.59 &  -  &  -  \\
 & 100/20  & 102.79 &  -  &  -  & 8.21 & 457 &  -  & 102.79 &  -  &  -  & 96.89 &  -  &  -  \\
 & 100/50  & 23.35 & 307 &  -  & 23.35 & 299 &  -  & 23.35 &  -  &  -  & 23.35 & 350 &  -  \\
 & 1000/01  & 25.02 & 0 &  -  & 28.63 & 14 &  -  & 25.02 & 0 &  -  & 25.02 & 0 &  -  \\
 & 1000/02  & 323.30 &  -  &  -  & 323.30 &  -  &  -  & 323.30 &  -  &  -  & 323.30 &  -  &  -  \\
 & 1000/05  & 533.37 &  -  &  -  & 533.37 &  -  &  -  & 533.37 &  -  &  -  & 533.37 &  -  &  -  \\
\midrule
 & 10/01  & 0.30 & 0 & 0 & 0.30 & 0 & 0 & 0.30 & 0 & 0 & 0.30 & 0 & 0 \\
 & 10/02  & 0.22 & 0 & 0 & 0.22 & 0 & 0 & 0.22 & 0 & 0 & 0.22 & 0 & 0 \\
 \parbox[t]{2mm}{\multirow{3}{*}{\rotatebox[origin=c]{90}{\texttt{IndicatorMode}}}}
 & 10/05  & 0.09 & 0 & 0 & 0.09 & 0 & 0 & 0.09 & 0 & 0 & 0.09 & 0 & 0 \\
 & 100/01  & 2.42 & 0 & 0 & 2.42 & 0 & 0 & 2.42 & 0 & 2 & 2.42 & 0 & 2 \\
 & 100/02  & 2.40 & 4 &  -  & 2.40 & 4 &  -  & 9.08 & 294 &  -  & 2.42 & 5 &  -  \\
 & 100/05  & 2.30 & 6 &  -  & 2.30 & 6 &  -  & 39.08 &  -  &  -  & 2.31 & 6 &  -  \\
& 100/10  & 79.59 &  -  &  -  & 79.59 &  -  &  -  & 79.59 &  -  &  -  & 79.59 &  -  &  -  \\
 & 100/20  & 102.79 &  -  &  -  &  -  &  -  &  -  & 102.79 &  -  &  -  & 102.79 &  -  &  -  \\
 & 100/50  & 23.35 &  -  &  -  & 23.35 &  -  &  -  & 23.35 &  -  &  -  & 23.35 & 576 &  -  \\
 & 1000/01  & 25.06 & 0 &  -  &  -  &  -  &  -  & 77.45 & 209 &  -  & 195.11 & 680 &  -  \\
 & 1000/02  & 323.30 &  -  &  -  &  -  &  -  &  -  & 323.30 &  -  &  -  &  -  &  -  &  -  \\
 & 1000/05  & 533.37 &  -  &  -  & 533.37 &  -  &  -  & 533.37 &  -  &  -  &  -  &  -  &  -  \\
\bottomrule
\end{tabular}
}
\caption{MIP solvers with \texttt{SOS1Mode} and \texttt{IndicatorMode}, Time in seconds (s), Gap in percent (\%).}
\label{table_sos1}
\end{table}

\begin{table}[!ht]
\resizebox{\textwidth}{!}{
\begin{tabular}{rr|rrr|rrr|rrr|rrr|rrr}
\toprule
 &  & \multicolumn{3}{c}{CPLEX} \vline & \multicolumn{3}{c}{Gurobi} \vline & \multicolumn{3}{c}{HiGHS} \vline & \multicolumn{3}{c}{SCIP} \vline & \multicolumn{3}{c}{Xpress}  \\
 & Inst & Obj & Gap  & Time  & Obj & Gap  & Time  & Obj & Gap  & Time  & Obj & Gap  & Time  & Obj & Gap  & Time  \\
\midrule
 & 10/01  & 0.30 & 0 & 0 & 0.30 & 0 & 0 & 0.30 & 0 & 0 & 0.30 & 0 & 0 & 0.30 & 0 & 0 \\
 & 10/02  & 0.22 & 0 & 0 & 0.22 & 0 & 0 & 0.22 & 0 & 0 & 0.22 & 0 & 0 & 0.22 & 0 & 0 \\
 \parbox[t]{2mm}{\multirow{3}{*}{\rotatebox[origin=c]{90}{\texttt{FortunyAmatMcCarlMode}}}}
 & 10/05  & 0.09 & 0 & 0 & 0.09 & 0 & 0 & 0.09 & 0 & 0 & 0.09 & 0 & 0 & 0.09 & 0 & 0 \\
 & 100/01  & 2.42 & 0 & 0 & 2.42 & 0 & 0 & 2.42 & 0 & 0 & 2.42 & 0 & 1 & 2.42 & 0 & 0 \\
 & 100/02  & 2.40 & 4 &  -  & 2.40 & 4 &  -  & 2.40 & 4 &  -  & 2.43 & 5 &  -  & 2.43 & 5 &  -  \\
 & 100/05  & 2.30 & 6 &  -  & 2.29 & 5 &  -  & 54.87 &  -  &  -  & 39.08 &  -  &  -  & 2.31 & 6 &  -  \\
 & 100/10  & 79.59 &  -  &  -  & 2.33 & 34 &  -  & 79.59 &  -  &  -  & 79.59 &  -  &  -  & 79.59 &  -  &  -  \\
 & 100/20  & 22.11 &  -  &  -  & 22.17 &  -  &  -  & 102.79 &  -  &  -  & 102.79 &  -  &  -  &  -  &  -  &  -  \\
 & 100/50  & 23.35 & 355 &  -  & 23.35 & 330 &  -  & 23.35 &  -  &  -  & 23.35 & 952 &  -  & 23.35 & 736 &  -  \\
 & 1000/01  & 25.02 & 0 &  -  & 25.02 & 0 &  -  & 25.02 & 0 &  -  & 70.13 & 180 &  -  & 25.02 & 0 &  -  \\
 & 1000/02  & 24.46 & 3 &  -  & 23.74 & 0 & 12 & 323.30 &  -  &  -  & 323.30 &  -  &  -  & 23.75 & 0 &  -  \\
 & 1000/05  & 533.37 &  -  &  -  & 533.37 &  -  &  -  & 533.37 &  -  &  -  & 533.37 &  -  &  -  & 533.37 &  -  &  -  \\
  \midrule
 & 10/01  & 0.30 & 0 & 344 & 0.30 & 0 & 17 & 0.30 & 0 & 269 & 0.30 & 0 & 83 & 0.30 & 0 & 447 \\
 & 10/02  & 1.82 & 739 &  -  & 0.22 & 0 &  -  & 3.12 &  -  &  -  & 0.33 & 0 & 495 & 0.30 & 38 &  -  \\
 \parbox[t]{2mm}{\multirow{3}{*}{\rotatebox[origin=c]{90}{\texttt{ProductMode}}}}
 & 10/05  & 0.53 &  -  &  -  & 7.22 &  -  &  -  & 8.72 &  -  &  -  & 0.67 &  -  &  -  & 0.60 &  -  &  -  \\
 & 100/01  & 18.08 & 647 &  -  & 14.37 & 494 &  -  & 23.28 & 0 & 52 & 21.87 & 803 &  -  & 20.74 & 0 & 55 \\
 & 100/02  & 27.06 &  -  &  -  &  -  &  -  &  -  &  -  &  -  &  -  &  -  &  -  &  -  &  -  &  -  &  -  \\
 & 100/05  &  -  &  -  &  -  & 48.56 &  -  &  -  &  -  &  -  &  -  &  -  &  -  &  -  &  -  &  -  &  -  \\
 & 100/10  &  -  &  -  &  -  & 75.06 &  -  &  -  &  -  &  -  &  -  &  -  &  -  &  -  & 78.62 &  -  &  -  \\
 & 100/20  &  -  &  -  &  -  & 99.82 &  -  &  -  &  -  &  -  &  -  &  -  &  -  &  -  & 101.51 &  -  &  -  \\
 & 100/50  & 247784.27 &  -  &  -  & 183.62 &  -  &  -  &  -  &  -  &  -  &  -  &  -  &  -  &  -  &  -  &  -  \\
 & 1000/01  & 45.15 & 80 &  -  &  -  &  -  &  -  &  -  &  -  &  -  & 58.77 & 135 &  -  & 147.68 & 490 &  -  \\
 & 1000/02  &  -  &  -  &  -  &  -  &  -  &  -  &  -  &  -  &  -  &  -  &  -  &  -  & 165.04 & 595 &  -  \\
 & 1000/05  &  -  &  -  &  -  &  -  &  -  &  -  &  -  &  -  &  -  &  -  &  -  &  -  &  -  &  -  &  -  \\
 \midrule
 & 10/01  & 0.30 & 0 & 512 & 0.30 & 0 &  -  & 0.30 & 0 &  -  & 0.30 & 0 &  -  & 0.30 & 0 &  -  \\
 & 10/02  & 0.22 & 0 & 0 & 0.22 & 0 & 0 & 0.22 & 0 & 7 & 0.22 & 0 & 36 & 0.22 & 0 & 52 \\
 \parbox[t]{2mm}{\multirow{3}{*}{\rotatebox[origin=c]{90}{\texttt{StrongDualityMode}}}}
 & 10/05  & 0.00 & 0 & 2 & 0.00 & 0 & 3 & 0.00 & 0 & 165 & 0.00 & 0 & 37 & 0.00 &  -  & 4 \\
 & 100/01  & 2.42 & 0 & 2 & 2.42 & 0 & 1 & 2.42 & 0 & 2 &  -  &  -  &  -  & 2.42 & 0 & 0 \\
 & 100/02  & 2.30 & 0 & 7 & 2.30 & 0 & 10 & 2.30 & 0 & 63 & 2.36 & 2 &  -  &  -  &  -  &  -  \\
 & 100/05  & 2.16 & 0 & 139 & 2.16 & 0 & 66 & 41.77 &  -  &  -  &  -  &  -  &  -  & 2.18 & 0 &  -  \\
 & 100/10  & 1.73 & 0 & 269 & 1.73 & 0 & 34 & 75.84 &  -  &  -  & 78.61 &  -  &  -  &  -  &  -  &  -  \\
 & 100/20  & 90.35 &  -  &  -  & 1.45 & 0 & 269 &  -  &  -  &  -  &  -  &  -  &  -  & 1.57 & 8 &  -  \\
 & 100/50  & 176.21 &  -  &  -  & 160.70 &  -  &  -  &  -  &  -  &  -  &  -  &  -  &  -  &  -  &  -  &  -  \\
 & 1000/01  & 25.01 & 0 & 153 & 25.01 & 0 & 22 & 25.01 & 0 & 32 & 25.01 & 0 & 20 & 25.01 & 0 & 10 \\
 & 1000/02  & 23.74 & 0 & 39 & 23.74 & 0 & 23 & 23.75 & 0 &  -  &  -  &  -  &  -  & 23.74 & 0 & 20 \\
 & 1000/05  & 529.71 &  -  &  -  & 24.43 & 0 & 543 & 408.46 &  -  &  -  &  -  &  -  &  -  & 24.43 & 0 & 482 \\
\bottomrule
\end{tabular}
}
\caption{MIP solvers with \texttt{FortunyAmatMcCarlMode}, \texttt{ProductMode} and \texttt{StrongDualityMode}, Time in seconds (s), Gap in percent (\%).}
\label{table_fa100}
\end{table}

\begin{table}[!ht]
\centering
\resizebox{11cm}{!}{
\begin{tabular}{rr|rrr|rrr|rrr}
\toprule
 &  & \multicolumn{3}{c}{Gurobi NonConvex} \vline & \multicolumn{3}{c}{Ipopt} \vline & \multicolumn{3}{c}{Knitro}  \\
 & Inst & Obj & Gap  & Time  & Obj & Gap  & Time  & Obj & Gap  & Time  \\
\midrule
 & 10/01  & 0.31 & 2 &  -  & 0.32 &  -  & 1 & 0.32 &  -  & 0 \\
 & 10/02  & 0.22 & 3 &  -  & 0.22 &  -  & 1 & 0.22 &  -  & 0 \\
\parbox[t]{2mm}{\multirow{3}{*}{\rotatebox[origin=c]{90}{\texttt{ProductMode}}}} & 10/05  & 0.67 &  -  &  -  & 0.09 &  -  & 0 & 0.09 &  -  & 0 \\
 & 100/01  & 2.42 & 0 & 6 & 2.43 &  -  & 18 & 2.46 &  -  & 0 \\
 & 100/02  & 2.71 & 17 &  -  & 2.41 &  -  & 12 & 2.43 &  -  & 0 \\
 & 100/05  & 54.87 &  -  &  -  & 2.47 &  -  & 24 & 2.54 &  -  & 0 \\
 & 100/10  & 79.59 &  -  &  -  & 2.64 &  -  & 12 & 2.35 &  -  & 0 \\
 & 100/20  & 102.79 &  -  &  -  & 3.43 &  -  & 11 & 3.51 &  -  & 0 \\
 & 100/50  & 185.45 &  -  &  -  & 19.46 &  -  & 22 & 19.46 &  -  & 0 \\
 & 1000/01  & 25.21 & 0 &  -  & 25.13 &  -  & 272 & 25.10 &  -  & 1 \\
 & 1000/02  & 323.30 &  -  &  -  & 24.08 &  -  & 148 & 23.84 &  -  & 2 \\
 & 1000/05  & 533.37 &  -  &  -  & 25.70 &  -  & 121 & 24.89 &  -  & 4 \\
 \midrule
 & 10/01  & 0.30 & 1 &  -  & 0.32 &  -  & 0 & 0.33 &  -  & 0 \\
 & 10/02  & 0.25 & 17 &  -  & 0.22 &  -  & 0 & 0.22 &  -  & 0 \\
\parbox[t]{2mm}{\multirow{3}{*}{\rotatebox[origin=c]{90}{\texttt{StrongDualityMode}}}} & 10/05  & 0.13 &  -  &  -  & 0.09 &  -  & 0 & 0.09 &  -  & 0 \\
 & 100/01  & 2.42 & 0 & 0 & 2.44 &  -  & 2 & 2.43 &  -  & 0 \\
 & 100/02  & 2.41 & 4 &  -  & 2.44 &  -  & 2 & 2.53 &  -  & 0 \\
 & 100/05  & 2.33 & 7 &  -  & 2.50 &  -  & 0 & 2.32 &  -  & 0 \\
 & 100/10  & 2.41 & 39 &  -  & 2.21 &  -  & 1 & 2.08 &  -  & 0 \\
 & 100/20  & 3.55 & 145 &  -  & 3.43 &  -  & 8 & 2.90 &  -  & 0 \\
 & 100/50  & 64.20 &  -  &  -  & 23.34 &  -  & 1 & 185.45 &  -  & 0 \\
 & 1000/01  & 25.03 & 0 &  -  & 54.89 &  -  & 3 & 25.08 &  -  & 3 \\
 & 1000/02  & 23.80 & 0 &  -  & 71.32 &  -  & 9 & 23.78 &  -  & 6 \\
 & 1000/05  & 24.86 & 1 &  -  & 29.52 &  -  & 175 & 24.75 &  -  & 4 \\
\bottomrule
\end{tabular}
}
\caption{Gurobi NonConvex and NLP solvers with \texttt{ProductMode} and \texttt{StrongDualityMode}, Time in seconds (s), Gap in percent (\%).}
\label{table_prod}
\end{table}

\clearpage

\ACKNOWLEDGMENT{%
The authors were partially supported by the Coordenação de Aperfeiçoamento de
Pessoal de Nível Superior - Brasil (CAPES) - Finance Code 001. The work of Alexandre Street was also partially supported by FAPERJ and CNPq.
}


\bibliographystyle{informs2014} 
\bibliography{ref.bib} 


\end{document}